\documentclass[]{article}
\begin{document}
\newtheorem{proposition}{Proposition}[section]
\newtheorem{definition}{Definition}[section]
\newtheorem{lemma}{Lemma}[section]

\title{\bf  Type $\theta$ Stokes' Theorem for Chains}
\author{Keqin Liu\\Department of Mathematics\\The University of British Columbia\\Vancouver, BC\\
Canada, V6T 1Z2}
\date{October, 2018}
\maketitle

\begin{abstract} We present type $\theta$ Stokes' theorem for type $\theta$ $k$-chains which extends the fundamental theorem of calculus in higher dimensions.
\end{abstract}

\medskip
The fundamental concepts and the important theorems in \cite{S} have  natural extentions if  the real number field $\mathcal{R}$ is replaced by the dual real number algebra $\mathcal{R}^{(2)}$ and the  real linear transformations  are replaced by $\mathcal{R}^{(2)}$-module maps. In this paper, we present type $\theta$ Stokes' theorem for type $\theta$ $k$-chains which extends the ordinary Stokes' theorem for chains.

\medskip
Section 1 gives the basic definitions about dual real $n$-space. Section 2 introduces dual real differentiation for dual real-valued functions with dual real multivariable by using  $\mathcal{R}^{(2)}$-module maps. Section 3 initiates the study of two types of integrations based on the two order relations on the dual real number algebra $\mathcal{R}^{(2)}$. Section 4 and Section 5 outline the notations and fundamental properties of alternating dual real tensors and dual  real differential forms. Section 6, which is the last section of this paper, defines the type $\theta$ integral of a dual real differentiable $k$-form and presents type $\theta$ Stokes' theorem for type $\theta$ $k$-chains.

\medskip
\section{Dual Real $n$-Space}

\medskip
Let $\mathcal{R}$ be the real number field with the identity $1$. Th dual real number algebra is the 2-dimensional real associative algebra
$\mathcal{R}^{(2)}=\mathcal{R}\oplus \mathcal{R}\,1^{\#}$ with the following product
$$
(x_1+x_2\,1^{\#})(y_1+y_2\,1^{\#})=x_1y_1+(x_1y_2+x_2y_1)\,1^{\#},
$$
where $x_1$, $x_2$, $y_1$, $y_2$ are real numbers, and $\{1, \,1^{\#}\}$ is a basis of the 2-dimensional real vector space $\mathcal{R}^{(2)}$. If
$x=x_1+x_2\,1^{\#}\in \mathcal{R}^{(2)}$ with $x_1$, $x_2\in \mathcal{R}$,
then $Re\,x:=x_1$ and $Ze\,x:=x_2$ are called the  real part and the  zero-divisor part of $x$, respectively. The dual real number algebra $\mathcal{R}^{(2)}$ is a commutative real associative algebras with zero-divisors.

\medskip
There are two generalized order relations on $\mathcal{R}^{(2)}$  which are compatible with the multiplication in $\mathcal{R}^{(2)}$.

\begin{definition}\label{def1.2} Let $x$ and $y$ be two elements of $\mathcal{R}^{(2)}$.
\begin{description}
\item[(i)] We say that $x$ is {\bf type 1 greater than} $y$ ( or $y$ is {\bf type 1 less than} $x$) and we write $x\stackrel{1}{>} y$ (or $y\stackrel{1}{<} x$) if
$$
\mbox{either}\quad \left\{\begin{array}{c}
Re\,x>Re\,y\\ Ze\,x\ge Ze\,y\end{array}\right.\quad\mbox{or}\quad
\left\{\begin{array}{c}
Re\,x=Re\,y\\ Ze\,x> Ze\,y\end{array}\right.
$$
\item[(ii)] We say that $x$ is {\bf type 2 greater than} $y$ ( or $y$ is {\bf type 2 less than} $x$) and we write $x\stackrel{2}{>} y$ (or $y\stackrel{2}{<} x$) if
$$
\mbox{either}\quad \left\{\begin{array}{c}
Re\,x>Re\,y\\ Ze\,x\le Ze\,y\end{array}\right.\quad\mbox{or}\quad
\left\{\begin{array}{c}
Re\,x=Re\,y\\ Ze\,x< Ze\,y\end{array}\right.
$$
\end{description}
\end{definition}

In the following of this paper, the letter $\theta$ always means an element in the set
$\{1, 2\}$. We use $x\stackrel{\theta}{\ge} y$ when $x\stackrel{\theta}{>} y$ or $x=y$.

\medskip
For any positive integer $n$, the {\bf dual real $n$-space} is the set
$$\mathcal{R}^{(2)n}:=
\{(x^1,  x^2,  \dots ,  x^n)\,|\, x^i\in \mathcal{R}^{(2)},\quad i=1, 2, \dots , n\}.$$
An element $x=(x^1,  x^2,  \dots ,  x^n)$ of $\mathcal{R}^{(2)n}$ is called a {\bf dual real vector} or a {\bf point}, and the dual real number $x^i=x^{i1}+x^{i2}\,1^{\#}$ with $x^{i1}$, $x^{i2}\in \mathcal{R}$ is called the {\bf $i$-th component} of $x$ for $1\le i\le n$. the dual real $n$-space $\mathcal{R}^{(2)n}$ is a left $\mathcal{R}^{(2)}$-module with respect to the following
addition and scalar multiplication:
$$
x+y:=(x^1+y^1,  x^2+y^2,  \dots ,  x^n+y^n), \quad ax:=(ax^1,  ax^2,  \dots ,  ax^n),
$$
where $x=(x^1,  x^2,  \dots ,  x^n)$, $y=(y^1,  y^2,  \dots ,  y^n)\in \mathcal{R}^{(2)n}$ and
$a\in \mathcal{R}^{(2)}$.

\medskip
Let $S$ be a non-empty subset  of left $\mathcal{R}^{(2)}$ -module $V$. An element of the left $\mathcal{R}^{(2)}$ -module $V$ is also called a dual real vector. $S$ is called
{\bf $\mathcal{R}^{(2)}$-linearly dependent} if there exist a finite number of distinct dual real vectors $v_1$, $v_2$, $\dots$, $v_n$ in $S$ and deal real numbers $c_1$, $c_2$, $\dots$, $c_n$, not all zero, such that $c_1v_1+c_2v_2+\dots+c_nv_n=0$. $S$ is called
{\bf $\mathcal{R}^{(2)}$-linearly independent} if $S$ is not $\mathcal{R}^{(2)}$-linearly dependent.  $S$ is called {\bf spanning set} of $v$ if every dual real vector in
$V$ is a {\bf $\mathcal{R}^{(2)}$-linearly combination} of  dual real vectors of $S$, i.e., for every
$v\in V$, there exist a finite number of distinct dual real vectors $v_1$, $v_2$, $\dots$, $v_n$ in $S$ and deal real numbers $c_1$, $c_2$, $\dots$, $c_n$ such that $v=c_1v_1+c_2v_2+\dots+c_nv_n$. $S$ is called a {\bf $\mathcal{R}^{(2)}$-basis} for the left $\mathcal{R}^{(2)}$ -module $V$ if $S$ is a $\mathcal{R}^{(2)}$-linearly independent spanning set of $V$.

\medskip
If $V$ is a free left $\mathcal{R}^{(2)}$-module having a finite $\mathcal{R}^{(2)}$-basis, then every
$\mathcal{R}^{(2)}$-basis for $V$ contains the same  number of dual real vectors. The unique number of  dual real vectors in each $\mathcal{R}^{(2)}$-basis for $V$ is called the
$\mathcal{R}^{(2)}$-{\bf dimension} of $V$ and is denoted by $\dim_{\mathcal{R}^{(2)}}(V)$.

\medskip
The dual real $n$-space $\mathcal{R}^{(2)n}$ is a free $\mathcal{R}^{(2)}$-module. Let $e_1:=(1, 0, 0, \dots ,  0)$, $e_2:=(0, 1, 0,  \dots ,  0)$, $\dots$,
$e_n:=(0, 0,  \dots ,  0, 1)$. Then $\{ e_1, e_2, \dots , e_n\}$ is a $\mathcal{R}^{(2)}$-basis for $\mathcal{R}^{(2)n}$ and is called the {\bf standard $\mathcal{R}^{(2)}$-basis}.

\bigskip
\section{Dual Real Differentiation}

\medskip
For $x=\displaystyle\sum_{i=1}^n (x^{i1}+1^{\#}x^{i2})e_i$ with $x^{i1}$, $x^{i2}\in \mathcal{R}$ for $1\le i\le n$, we define
\begin{equation}\label{eq1.3}
||x||:=\Big(2\displaystyle\sum_{i=1}^n (x^{i1})^2+
\displaystyle\sum_{i=1}^n (x^{i2})^2\Big)^{\frac12}.
\end{equation}
Then  the real-valued function $||\,\,||: \mathcal{R}^{(2)n}\to \mathcal{R}$ defined by (\ref{eq1.3}) is a norm on the dual real $n$-space $\mathcal{R}^{(2)n}$.

\medskip
The dual real $n$-space $\mathcal{R}^{(2)n}$ is a metric space with the distance function $||\,\,||$ defined by (\ref{eq1.3}). If $a\in\mathcal{R}^{(2)n}$ and $\epsilon\in \mathcal{R}$, we use $N(a;\,\epsilon)$  and $N^*(a;\,\epsilon)$ to denote the ordinary $\epsilon$-{\bf neighborhood} and {\bf deleted $\epsilon$-neighborhood} of $a$, respectively, i.e.,
$$N(a;\,\epsilon)=\{x\in \mathcal{R}^{(2)n}\,|\, ||x-a||<\epsilon\}\quad\mbox{and}\quad
N^*(a;\,\epsilon):=N(a;\,\epsilon)\setminus\{a\}.$$

\medskip
\begin{definition}\label{def2.1} Let $U$ be an open subset of $\mathcal{R}^{(2)n}$. A function $f: U\to \mathcal{R}^{(2)m}$ is
{\bf dual real differentiable} at $a\in U$ if there is a
$\mathcal{R}^{(2)}$-module map $\lambda: \mathcal{R}^{(2)n}\to \mathcal{R}^{(2)m}$ such that
$$
\lim_{x\to a}\frac{||f(x)-f(a)-\lambda (x-a)||}{||x-a||}=0,
$$
i.e., the following property holds:
\begin{eqnarray*}
&&\mbox{for each positive real number $\epsilon>0$, there exist a  positive real} \nonumber\\
&&\mbox{number $\delta>0$  such that}\nonumber\\
&& x\in N^*(a;\,\delta)\Rightarrow
\frac{||f(x)-f(a)-\lambda (x-a)||}{||x-a||}<\epsilon.
\end{eqnarray*}
The $\mathcal{R}^{(2)}$-module map $\lambda$ is denoted by $\mathcal{D}f(a)$ and
called the {\bf dual real derivative} of $f$ at $a$. If $f$ is dual real differentiable at each point of the open subset $U$, then $f$ is said to be {\bf  dual real differentiable} on $U$.
\end{definition}

\medskip
\begin{definition}\label{def2.2} Let $D$ be a subset of $\mathcal{R}^{(2)n}$. A function
$f: D\to \mathcal{R}^{(2)m}$ is said to be
{\bf continuous} at $a\in D$ if for each positive real number $\epsilon>0$, there exists a  positive real number $\delta>0$  such that
$$x\in N(a;\,\delta)\cap D\Rightarrow ||f(x)-f(a)||<\epsilon .$$
If $f$ is continuous at each point of the open subset $D$, then $f$ is said to be
{\bf continuous} on $D$.
\end{definition}

\medskip
One can check that if $f:\mathcal{R}^{(2)n}\to \mathcal{R}^{(2)m}$ is  dual real differentiable at $a\in \mathcal{R}^{(2)n}$, then $f$ is continuous at $a$. Also, the composition of two dual real differentiable functions is  dual real differentiable. In fact, we have

\medskip
\begin{proposition}\label{pr2.3} {\bf (Chain Rule)} If
$f: \mathcal{R}^{(2)n}\to \mathcal{R}^{(2)m}$ is dual real differentiable at $a\in\mathcal{R}^{(2)n}$, and $g: \mathcal{R}^{(2)m}\to \mathcal{R}^{(2)p}$ is dual real differentiable at $f(a)\in\mathcal{R}^{(2)m}$, then the composition
$g\circ f: \mathcal{R}^{(2)n}\to \mathcal{R}^{(2)p}$ is dual real differentiable at $a$, and
$$ \mathcal{D}(g\circ f)(a)=\mathcal{D}g \big(f(a)\big)\circ \mathcal{D}f(a). $$
\end{proposition}

\hfill\raisebox{1mm}{\framebox[2mm]{}}

\bigskip
\section{Type $\theta$ Integrations}

\medskip
Let $a^i$, $b^i\in \mathcal{R}^{(2)}$ with $a^i\stackrel{\theta}{<}b^i$ for $1\le i\le n$. Recall that the  type $\theta$ closed interval $[a^i, b^i]_{\theta}$ is defined by
$
[a^i, b^i]_{\theta}:=\big\{\, x\in \mathcal{R}^{(2)}\,|\, a^i\stackrel{\theta}{\le}x\stackrel{\theta}{\le}b^i\,\big\}.
$
The subset $[a^1, b^1]_{\theta}\times \cdots \times [a^n, b^n]_{\theta}$ of $\mathcal{R}^{(2)n}$
is called a {\bf type $\theta$ closed rectangle}. A {\bf partition} of a type $\theta$ closed rectangle  $[a^1, b^1]_{\theta}\times \cdots \times [a^n, b^n]_{\theta}$ is a collection
$P=(P_1, \dots , P_n)$, where $P_i$ is a partition of  $[a^i, b^i]_{\theta}$ with the partition points:
$$
a^i=x^{i(0)}\stackrel{\theta}{<} x^{i(1)}\stackrel{\theta}{<} \dots \stackrel{\theta}{<} x^{i(N_i)}=b^i\quad\mbox{for $1\le i\le n$.}
$$
The partition $P=(P_1, \dots , P_n)$ produces $N_1N_2\cdots N_n$ type $\theta$ closed subrectangles:
$$
S=S_{j_1, \dots , j_n}=[x^{1(j_1-1)}, x^{1(j_1)}]_{\theta}\times [x^{2(j_2-1)}, x^{2(j_2)}]_{\theta}\times\cdots \times [x^{n(j_n-1)}, x^{n(j_n)}]_{\theta} ,
$$
where $1\le j_1\le N_1$, $\dots$ , $1\le j_n\le N_n$, which are called the type $\theta$
{\bf subrectangles of the partition $P$}. The {\bf dual real volume}
$v(S)=v(S_{j_1, \dots , j_n})$ is defined by
$$
v(S)=v(S_{j_1, \dots , j_n})=\Big(x^{1(j_1)}-x^{1(j_1-1)}\Big)\cdots
\Big(x^{n(j_n)}-x^{n(j_n-1)}\Big).
$$
Clearly, $v(S)=v(S_{j_1, \dots , j_n})\stackrel{\theta}{\ge}0$.

\medskip
If $f: A\to \mathcal{R}^{(2)}$ is a function on a  type $\theta$ closed rectangle
$A$ in $\mathcal{R}^{(2)n}$, then
\begin{equation}\label{eq3.1}
f(x)=f_{Re}(x)+f_{Ze}(x)\,1^{\#}\quad\mbox{for $x\in A$},
\end{equation}
where $x=(x^1, \dots , x^n)=(x^{11}+x^{12}\,1^{\#}, \dots , x^{n1}+x^{n2}\,1^{\#})$ and
$$
f_{\clubsuit}(x)=f_{\clubsuit}(x^{11}, x^{12}, x^{21}, x^{22}, \dots , x^{n1}, x^{n2})
\quad\mbox{with $\clubsuit\in \{Re, \, Ze\}$}
$$
is a real-valued function of $2n$ real variables $x^{11}$, $x^{12}$, $\dots$, $x^{n1}$, $x^{n2}$.
We say that the function $f: A\to \mathcal{R}^{(2)}$ given by (\ref{eq3.1}) is {\bf bounded} if both $f_{Re}$ and $f_{Ze}$ are bounded on $A$ ($\subseteq \mathcal{R}^{(2)n}$).

\medskip
Suppose that $A$ is a type $\theta$ closed rectangle in $\mathcal{R}^{(2)n}$,
$f: A\to \mathcal{R}^{(2)n}$ is a bounded function, and $P=(P_1, \dots , P_n)$ is a partition of $A$. For each type $\theta$ subrectangle $S$ of the partition $P$, let
$f\big|S: S\to \mathcal{R}^{(2)}$ be the restriction of $f$ to $S$.
The range of $f\big|S$ is a bounded subset of the real number field $\mathcal{R}$. Hence, both
$$\sup_Sf_{\clubsuit}:=\sup \big\{f_{\clubsuit}(x)\,|\, x\in S\big\}$$
and
$$\inf_Sf_{\clubsuit}:=\inf \big\{f_{\clubsuit}(x)\,|\, x\in S\big\}$$
exist for $\clubsuit\in \{Re, \, Ze\}$ and each type $\theta$ subrectangle $S$ of the partition $P$. Since $f_{\clubsuit}$ is bounded on $A$, there exist  real numbers $m$, $M\in \mathcal{R}$ such that
$$
 m\le f_{\clubsuit}(x)\le M\quad\mbox{for all $x\in A$ and $\clubsuit\in\{Re, \, Ze\}$}.
$$
Clearly, we have
\begin{equation}\label{eq3.2}
 \sup_Sf_{\clubsuit}\le M\quad \mbox{and}\quad
\inf_Sf_{\clubsuit}\ge m
\end{equation}
for $\clubsuit\in\{Re, \, Ze\}$ and any type $\theta$ subrectangle $S$ of the partition $P$.

\medskip
We define the {\bf type $\theta$ upper sum }$U_{\theta}(P,  f)$ of $f$ with respect to the partition $P$ of  the type $\theta$ rectangle $A$ to be
\begin{eqnarray}\label{eq3.3}
U_{\theta}(P,  f)&=&\left\{\begin{array}{c}
\displaystyle\sum_{S\,\widehat{\in}\,P}\big(\sup_Sf_{Re}+1^{\#}\,\sup_Sf_{Ze}\big)\,v(S)
\quad\mbox{$\theta =1$;}\\  \\
\displaystyle\sum_{S\,\widehat{\in}\,P}\big(\sup_Sf_{Re}+1^{\#}\,\inf_Sf_{Ze}\big)\,v(S)
\quad\mbox{$\theta =2$}
\end{array}\right.
\end{eqnarray}
and the {\bf type $\theta$ lower sum} $L_{\theta}(P,  f)$ of $f$ with respect to the partition $P$ of  the type $\theta$ rectangle $A$ to be
\begin{eqnarray}\label{eq3.4}
L_{\theta}(P,  f)&=&\left\{\begin{array}{c}
\displaystyle\sum_{S\,\widehat{\in}\,P}\big(\inf_Sf_{Re}+1^{\#}\,\inf_Sf_{Ze}\big)\,v(S)
\quad\mbox{$\theta =1$;}\\  \\
\displaystyle\sum_{S\,\widehat{\in}\,P}\big(\inf_Sf_{Re}+1^{\#}\,\sup_Sf_{Ze}\big)\,v(S)
\quad\mbox{$\theta =2$},
\end{array}\right.
\end{eqnarray}
where the notation $S\,\widehat{\in}\,P$ means that $S$ is a  type $\theta$ subrectangle $S$ of the partition $P$.

\medskip
Using (\ref{eq3.3}) and (\ref{eq3.4}), we have
\begin{eqnarray}\label{eq3.13}
&&\clubsuit \Big(\big(m+1^{\#}m\big)\, v(A)\Big)\stackrel{1}{\le}
\clubsuit L_{1}(P,  f)\nonumber\\
&&\qquad \stackrel{1}{\le}
\clubsuit U_{1}(P,f)\stackrel{1}{\le}\clubsuit\Big(\big(M+1^{\#}M\big)\,v(A)\Big)
\end{eqnarray}
for $\clubsuit\in\{Re, \, Ze\}$ and
\begin{eqnarray}\label{eq3.14}
&&Re\Big(\big(m+1^{\#}M\big)\, v(A)\Big)\stackrel{2}{\le}Re\,L_{2}(P,  f)\nonumber\\
&&\qquad \stackrel{2}{\le}
Re\,U_{2}(P,f)\stackrel{2}{\le}Re\Big(\big(M+1^{\#}m\big)\,v(A)\Big)
\end{eqnarray}
and
\begin{eqnarray}\label{eq3.15}
&&Ze\Big(\big(m+1^{\#}M\big)\, v(A)\Big)\stackrel{2}{\ge}Ze\,L_{2}(P,  f)\nonumber\\
&&\qquad \stackrel{2}{\ge}
Ze\,U_{2}(P,f)\stackrel{2}{\ge}Ze\Big(\big(M+1^{\#}m\big)\,v(A)\Big).
\end{eqnarray}

Let $\mathcal{P}_{\theta}$ be the set of all partitions of the type $\theta$ closed
rectangle $A=[a^1, b^1]_{\theta}\times \cdots \times [a^n, b^n]_{\theta}$, i.e.,
$$
\mathcal{P}_{\theta}:=\Big\{P\,\Big|\, \mbox{$P$ is a partition of
$[a^1, b^1]_{\theta}\times \cdots \times [a^n, b^n]_{\theta}$}\Big\}.
$$
By (\ref{eq3.13}), (\ref{eq3.14}) and (\ref{eq3.15}), the following four sets
$$
\big\{\clubsuit U_{\theta}(P,  f)\,\big|\, P\in \mathcal{P}_{\theta}\big\},\quad
\big\{\clubsuit L_{\theta}(P,  f)\,\big|\, P\in \mathcal{P}_{\theta}\big\}\quad\mbox{with}\quad
\clubsuit\in\{Re, \, Ze\}
 $$
are bounded subsets of the real number field $\mathcal{R}$. Hence, the supremums and infimums of the four sets exist. We now define the {\bf type $\theta$ lower integral}
$$\underline{\displaystyle\int_A} f(x)d_{\theta}V
=\underline{\displaystyle\int_A} f(x^1, \dots , x^n)\,dx^1\cdots dx^n$$
and  the {\bf type $\theta$ upper integral}
$$\overline{\displaystyle\int_A} f(x)d_{\theta}V
=\overline{\displaystyle\int_A} f(x^1, \dots , x^n)\,dx^1\cdots dx^n$$
of $f(x)$ on the type $\theta$ closed
rectangle $A=[a^1, b^1]_{\theta}\times \cdots \times [a^n, b^n]_{\theta}$  by
\begin{eqnarray*}\label{eq3.16}
&&\underline{\displaystyle\int_A} f(x)d_{1}V\nonumber\\
&=&\sup\{Re L_{1}(P,  f)\,\big|\, P\in \mathcal{P}_{1}\}+
1^{\#}\sup\{Ze L_{1}(P,  f)\,\big|\, P\in \mathcal{P}_{1}\},
\end{eqnarray*}
\begin{eqnarray*}\label{eq3.17}
&&\underline{\displaystyle\int_A} f(x)d_{2}V\nonumber\\
&=&\sup\{Re L_{2}(P,  f)\,\big|\, P\in \mathcal{P}_{2}\}+
1^{\#}\inf\{Ze L_{2}(P,  f)\,\big|\, P\in \mathcal{P}_{2}\},
\end{eqnarray*}
\begin{eqnarray*}\label{eq3.18}
&&\overline{\displaystyle\int_A} f(x)d_{1}V\nonumber\\
&=&\inf\{Re U_{1}(P,  f)\,\big|\, P\in \mathcal{P}_{1}\}+
1^{\#}\inf\{Ze U_{1}(P,  f)\,\big|\, P\in \mathcal{P}_{1}\}
\end{eqnarray*}
and
\begin{eqnarray*}\label{eq3.19}
&&\overline{\displaystyle\int_A} f(x)d_{2}V\nonumber\\
&=&\inf\{Re U_{2}(P,  f)\,\big|\, P\in \mathcal{P}_{2}\}+
1^{\#}\sup\big\{Ze U_{2}(P,  f)\,\big|\, P\in \mathcal{P}_{2}\}.
\end{eqnarray*}

If the  type $\theta$ lower integral and the type $\theta$ upper integral of $f(x)$ on
$A=[a^1, b^1]_{\theta}\times \cdots \times [a^n, b^n]_{\theta}$ are equal, i.e.,
if
$$\underline{\displaystyle\int_A} f(x)d_{\theta}V=
\overline{\displaystyle\int_A} f(x)d_{\theta}V,$$
then we say that $f$ is {\bf type $\theta$ integrable} on $A$, we denote their common value by
$\displaystyle\int_A f(x)d_{\theta}V$ which is called the {\bf type $\theta$ integral} of $f$  on
$A$.

\medskip
A basic fact is that if $f: A\to \mathcal{R}^{(2)}$ is a continuous function on a  type $\theta$ closed rectangle $A$ in $\mathcal{R}^{(2)n}$, then  $f$ is  type $\theta$ integrable on $A$.

\medskip
\section{Alternating Dual Real Tensors}

Let $V$ be a left $\mathcal{R}^{(2)}$-module, and let $k$ be a positive integer. We define $V^k$ by $V^k:=\underbrace{V\times \cdots\times V}_{k}$. A function $T: V^k\to \mathcal{R}^{(2)}$ is called a {\bf dual real $k$-tensor} if $T$ is $\mathcal{R}^{(2)}$-{\bf multilinear}, i.e., if for each $i$ with $1\le i\le k$ we have
\begin{eqnarray*}
&&T(v_1, \dots , v_{i-1}, v_i+v_i',  v_{i+1}, \dots , v_k)\\
&=&T(v_1, \dots , v_{i-1}, v_i,  v_{i+1}, \dots , v_k)
+T(v_1, \dots , v_{i-1}, v_i',  v_{i+1}, \dots , v_k)
\end{eqnarray*}
and
$$
T(v_1, \dots , v_{i-1}, av_i,  v_{i+1}, \dots , v_k)
=a\,T(v_1, \dots , v_{i-1}, v_i,  v_{i+1}, \dots , v_k),
$$
where $v_i$, $v_i'\in V$ and $a\in \mathcal{R}^{(2)}$. Let
$$
\mathbf{T}^k(V):=\{ T\,|\, \mbox{$T: V^k\to \mathcal{R}^{(2)}$ is called a
 dual real $k$-tensor}\}.
$$
Then $\mathbf{T}^k(V)$ becomes a left $\mathcal{R}^{(2)}$-module if for $S$, $T\in\mathbf{T}^k(V)$ and $a\in \mathcal{R}^{(2)}$ we define
\begin{eqnarray*}
\label{eq4.1}(S+T)(v_1, \dots , v_k):&=&S(v_1, \dots , v_k)+T(v_1, \dots , v_k)\\
\label{eq4.2}(aS)(v_1, \dots , v_k):&=&a\cdot S(v_1, \dots , v_k).
\end{eqnarray*}

\medskip
The {\bf tensor product} $S\otimes T\in \mathbf{T}^{k+\ell}(V)$ of $S\in\mathbf{T}^k(V)$ and $T\in\mathbf{T}^{\ell}(V)$ is defined by
\begin{eqnarray*}\label{eq4.3}
&&(S\otimes T)(v_1, \dots , v_k, v_{k+1}, \dots ., v_{k+\ell}):\nonumber\\
&=&S(v_1, \dots , v_k)\cdot T(v_{k+1}, \dots ., v_{k+\ell}).
\end{eqnarray*}

\medskip
A dual real $k$-tensor $w\in \mathbf{T}^k(V)$ is called {\bf alternating} if
\begin{eqnarray*}\label{eq4.18}
&&w(v_1, \dots , v_{i-1}, v_i,  v_{i+1}, \dots , v_{j-1}, v_j,  v_{j+1}, \dots , v_k)\nonumber\\
&=&-w(v_1, \dots , v_{i-1}, v_j,  v_{i+1}, \dots , v_{j-1}, v_i,  v_{j+1}, \dots , v_k).
\end{eqnarray*}
The set of all alternating dual real $k$-tensors is a left $\mathcal{R}^{(2)}$-submodule of
$\mathbf{T}^k(V)$, which is denoted by $\mathbf{\Lambda}^k(V)$.

\medskip
If $T\in \mathbf{T}^k(V)$, we define $Alt(T)$ by
$$
Alt(T)(v_1, \dots , v_k):=\displaystyle\frac{1}{k!}\sum_{\sigma\in\mathbf{S}_k}
(sgn\,\sigma)\,T(v_{\sigma (1)}, \dots , v_{\sigma (k)}),
$$
where $\mathbf{S}_k$ is the set of all permutations on the set $\{1, 2, \dots , k\}$, and $sgn\,\sigma$ is the sign of a permutation $\sigma\in\mathbf{S}_k$.

\medskip
For $w\in\mathbf{\Lambda}^k(V)$ and $\eta\in\mathbf{\Lambda}^{\ell}(V)$, the {\bf wedge product}
$ w\wedge \eta\in \mathbf{\Lambda}^{k+\ell}(V)$ is defined by
$$
w\wedge \eta:=\displaystyle\frac{(k+\ell)!}{k!\,\ell !}\,Alt(w\otimes \eta).
$$

\medskip
We finish this section with the following

\medskip
\begin{proposition}\label{pr4.6} Let $v_1$, $\dots$, $v_n$ be a $\mathcal{R}^{(2)}$-basis for a
left $\mathcal{R}^{(2)}$-module $V$. If $\varphi _1$, $\dots$, $\varphi _n$ is the  dual
$\mathcal{R}^{(2)}$-basis of $V$ with respect to the $\{v_1, \dots, v_n\}$, i.e., $\varphi _i: V\to \mathcal{R}^{(2)}$ is a $\mathcal{R}^{(2)}$-module map satisfying $\varphi _i(v_j): =\delta_{ij}$ for
$1\le i, \, j\le n$, then the set of  all wedge products
$$
\Big\{\varphi _{i_1}\wedge \varphi _{i_2}\wedge\cdots \wedge \varphi _{i_k}\,\Big|\,
1\le i_1<i_2<\dots < i_k\le n\Big\}
$$
is a $\mathcal{R}^{(2)}$-basis of $\mathbf{\Lambda}^k(V)$. Therefore,
$\dim_{\mathcal{R}^{(2)}}\big(\mathbf{\Lambda}^k(V)\big)=
\left(\begin{array}{c}n\\k\end{array}\right)=\displaystyle\frac{n!}{k!(n-k)!}.$
\end{proposition}

\medskip
\section{Dual  Real Differential Forms}

For $p\in \mathcal{R}^{(2)n}$, let $\mathcal{R}^{(2)n}_p$ be the set defined by
$\mathcal{R}^{(2)n}_p:=\{\, (p, v)\,|\, v\in \mathcal{R}^{(2)n}\}$. The pair  $(p, v)$ is also denoted by $v_p$. $\mathcal{R}^{(2)n}_p$ is made into a left $\mathcal{R}^{(2)}$-module  by defining
\begin{eqnarray*}
\label{eq4.54} (p, v)+(p, w):&=&(p, v+w),\\
\label{eq4.55} a\cdot (p, v):&=&(p, av)
\end{eqnarray*}
for $v$, $w\in  \mathcal{R}^{(2)n}$ and $a\in  \mathcal{R}^{(2)}$. Let $\{e_1, e_2, \dots,
e_n\}$ be the standard $\mathcal{R}^{(2)}$-basis for $\mathcal{R}^{(2)n}$. Then
$\{(e_1)_p, (e_2)_p, \dots, (e_n)_p\}$ is a $\mathcal{R}^{(2)}$-basis for $\mathcal{R}^{(2)n}_p$, which is called the {\bf standard $\mathcal{R}^{(2)}$-basis} for  $\mathcal{R}^{(2)n}_p$.

\medskip
A function $w: \mathcal{R}^{(2)n}\to
\displaystyle\bigcup_{q\in \mathcal{R}^{(2)n}} \mathbf{\Lambda}^k(\mathcal{R}^{(2)n}_q)$ is called a {\bf dual real $k$-differential form} if there exist dual real differentiable functions
$w_{i_1, \dots, i_k}: \mathcal{R}^{(2)n}\to\mathcal{R}^{(2)}$ for $1\le i_1<\dots <i_k\le n$ such that
\begin{equation}\label{eq4.56}
w(p)=\displaystyle\sum_{1\le i_1< \dots < i_k\le n}w_{i_1, \dots, i_k}(p)\cdot
\Big(\varphi _{i_1}(p)\wedge \cdots \wedge \varphi _{i_k}(p)\Big)\in
\mathbf{\Lambda}^k(\mathcal{R}^{(2)n}_p)
\end{equation}
for all $p\in \mathcal{R}^{(2)n}$, where $\{\varphi _{i_1}(p), \dots, \varphi _{i_n}(p)\}$ is the dual $\mathcal{R}^{(2)}$-basis with respect to the standard $\mathcal{R}^{(2)}$-basis
$\{(e_1)_p, (e_2)_p, \dots, (e_n)_p\}$ of  $\mathcal{R}^{(2)n}_p$.

\medskip
If $f: \mathcal{R}^{(2)n}\to \mathcal{R}^{(2)}$ is dual real  differentiable, then
$\mathcal{D}f(p)\in \mathbf{\Lambda}^1(\mathcal{R}^{(2)n})$ for all $p\in \mathcal{R}^{(2)n}$, and we can define a dual real $1$-form $df$ by
\begin{equation}\label{eq4.60}
(df)(p)(v_p):=\mathcal{D}f(p)(v)\quad\mbox{for $v_p\in \mathcal{R}^{(2)n}_p$}.
\end{equation}
By (\ref{eq4.60}), we know that $(df)(p)\in \mathbf{\Lambda}^1(\mathcal{R}^{(2)n}_p)$ for
$p\in\mathcal{R}^{(2)n}$.

\medskip
Since the $i$th-projection function
$\pi^i:\mathcal{R}^{(2)n}\to \mathcal{R}^{(2)}$ is $\mathcal{R}^{(2)}$-linear, we get
$$
\mathcal{D}\pi^i(p)=\pi^i\quad\mbox{for $p\in\mathcal{R}^{(2)n}$}.
$$
After denoting the dual real differential
$1$-form $d\pi^i$ by $dx^i$, we have
\begin{equation}\label{eq4.62}
(d x^i)(p)(v_p)=(d\pi^i)(p)(v_p)=(D\pi^i)(p)(v)=\pi^i(v)=v^i,
\end{equation}
where $v=(v^1, \dots, v^n)$ and $v^i\in \mathcal{R}^{(2)}$ for $1\le i\le n$. It follows from
(\ref{eq4.62}) that $(dx^i)(p)=\varphi _{i}(p)$, i.e., $\{dx^1(p), \dots, dx^n(p)\}$ is just
the dual real $\mathcal{R}^{(2)}$-basis  of the  standard $\mathcal{R}^{(2)}$-basis
$\{(e_1)_p, (e_2)_p, \dots, (e_n)_p\}$ for  $\mathcal{R}^{(2)n}_p$. Hence, every dual real
differential $k$-form $w$ given by (\ref{eq4.56}) can be written as
$$
w=\displaystyle\sum_{1\le i_1< \dots < i_k\le n}w_{i_1, \dots, i_k}\,
d x^{i_1}\wedge d x^{i_2}\wedge\cdots \wedge d x^{i_k}.
$$
The {\bf dual real differential $dw$} of $w$ is a dual real differentiable $(k+1)$-form on $\mathcal{R}^{(2)n}$, where $dw$ is defined by
$$
dw=\displaystyle\sum_{1\le i_1<\dots<i_k\le n}(dw_{i_1, \dots, i_k})\,\wedge dx^{i_1}\wedge \dots dx^{i_k}.
$$

\medskip
If $f: \mathcal{R}^{(2)n}\to \mathcal{R}^{(2)m}$ is dual real  differentiable, we get a $\mathcal{R}^{(2)}$-linear map
$f^*: \mathbf{\Lambda}^k(\mathcal{R}^{(2)m}_{f(p)})\to \mathbf{\Lambda}^k(\mathcal{R}^{(2)n}_{p})$, where $f^*$ is given by
$$
(f^*T)\Big((v_1)_p, \dots, (v_k)_p\Big)
=T\Big(\big((\mathcal{D}f)(p)(v_1)\big)_{f(p)}, \dots,
\big((\mathcal{D}f)(p)(v_k)\big)_{f(p)}\Big),
$$
where $T\in \mathbf{\Lambda}^k(\mathcal{R}^{(2)m}_{f(p)})$ and $v_1$, $\dots$,
$v_k\in \mathcal{R}^{(2)n}$.
We can therefore define a
dual real  differentiable $k$-form on $\mathcal{R}^{(2)n}$ by
$$
(f^*w)(p):=f^*\Big(w\big(f(p)\big)\Big)\quad\mbox{for $p\in\mathcal{R}^{(2)n}$},
$$
which means that if $v_1$, $\dots$,  $v_k\in \mathcal{R}^{(2)n}$, then
$(f^*w)(p)\in \mathbf{\Lambda}^k(\mathcal{R}^{(2)n}_p)$ and
\begin{eqnarray}\label{eq4.70}
&&(f^*w)(p)\Big((v_1)_p, \dots, (v_k)_p\Big)\nonumber\\
&=&w\big(f(p)\Big(\big((\mathcal{D}f)(p)(v_1)\big)_{f(p)}, \dots, \big((\mathcal{D}f)(p)(v_k)\big)_{f(p)}\Big).
\end{eqnarray}

\medskip
\section{Type $\theta$ Stokes' Theorem for Type $\theta$ $k$-Chain}

For a non-negative integer $k$, we define $[0,\, b]_{\theta}^0:=\{0\}$ and
$$
[0,\, b]_{\theta}^k:=\underbrace{[0,\,b]_{\theta}\times \cdots \times [0,\,b]_{\theta}}_k
\quad\mbox{for $k\ge 1$},
$$
where $\theta\in\{1, 2\}$,  $[0,\,b]_{\theta}\subseteq \mathcal{R}^{(2)}$ is a type $\theta$ closed interval, $b=1+(-1)^{\theta -1}r1^{\#}$ and $r\ge 0$ is a nonnegative real number. A
{\bf type $\theta$ singular $k$-cube} in
$A\subseteq \mathcal{R}^{(2)n}$ is a function $c^{b, \theta}: [0,\, b]_{\theta}^k\to A$ which extends to be a dual real differentiable function on some open set $U$ of $\mathcal{R}^{(2)k}$ such that $[0,\, b]_{\theta}^k\subseteq U$.
A finite sum of type $\theta$  singular $k$-cubes in $A$ with integer coefficients is called a {\bf type $\theta$ $k$-chain} in $A$.
A type $\theta$ singular $0$-cube in $A$ is a point in $A\subseteq \mathcal{R}^{(2)n}$. The
{\bf type $\theta$ standard singular $k$-cube} is the function
$I^{k, b, \theta}: [0,\, b]_{\theta}^k \to \mathcal{R}^{(2)k}$ defined by
$$
I^{k, b, \theta}(x):=x \quad\mbox{for $x\in [0,\, b]_{\theta}^k$}.
$$

\medskip
For $k\ge 1$, $j\in\{1, \dots , k\}$ and $\alpha\in\{0, \, 1\}$, we define a type $\theta$ singular $(k-1)$-cube $I^{k, b, \theta}_{(j, \alpha)}:  [0,\, b]_{\theta}^{k-1}\to
[0,\, b]_{\theta}^k\subseteq \mathcal{R}^{(2)k}$ by
\begin{eqnarray}\label{eq4.111}
I^{k, b, \theta}_{(j, \alpha)}(x):&=& I^{k, b, \theta}(x^1, \dots , x^{j-1},
\stackrel{\stackrel{j-th}{\downarrow}}{\alpha \,b}, x^j, \dots , x^{k-1})\nonumber\\
&=&(x^1, \dots , x^{j-1},
\stackrel{\stackrel{j-th}{\downarrow}}{\alpha \,b}, x^j, \dots , x^{k-1})
\in [0,\, b]_{\theta}^k,
\end{eqnarray}
where $x=(x^1, \dots , x^{k-1})\in [0,\, b]_{\theta}^{k-1}$ and $x^1$, $\dots$,
$x^{k-1}\in  [0,\, b]_{\theta}$. $I^{k, b, \theta}_{(j, 0)}$ and $I^{k, b, \theta}_{(j, 1)}$ are called the {\bf
$(j, 0)$-face} of $I^{k, b, \theta}$ and {\bf $(j, 1)$-face} of $I^{k, b, \theta}$, respectively.

\medskip
The {\bf boundary} $\partial I^{k, b, \theta}$ of $I^{k, b, \theta}$ is the  type $\theta$
$(k-1)$-chain defined by
$$
\partial I^{k, b, \theta}:=\displaystyle\sum_{j=1}^k\displaystyle\sum_{\alpha=0, 1}
(-1)^{j+\alpha}\, I^{k, b, \theta}_{(j, \alpha)}.
$$
For a general type $\theta$ singular $k$-cube
$c^{b, \theta}: [0,\, b]_{\theta}^k\to A\subseteq \mathcal{R}^{(2)n}$, the
{\bf $(j, \alpha)$-face} $c^{b, \theta}_{j, \alpha}$ of $c^{b, \theta}$ is defined by
$$
c^{b, \theta}_{j, \alpha}:=c^{b, \theta}\circ \Big(I^{k, b, \theta}_{(j, \alpha)}\Big):
[0,\, b]_{\theta}^{k-1}\to A\subseteq \mathcal{R}^{(2)n}
$$
and the {\bf boundary} $\partial c^{b, \theta}$ of $c^{b, \theta}$ is defined by
$$
\partial c^{b, \theta}:=\displaystyle\sum_{j=1}^k\displaystyle\sum_{\alpha=0, 1}
(-1)^{j+\alpha}\, c^{k, \theta}_{(j, \alpha)}.
$$
If $a_i$ is an integer and $c^{b, \theta}_i$ is a  type $\theta$ singular $k$-cube for
$1\le i\le m$, then the {\bf boundary} of type $\theta$
$k$-chain $\displaystyle\sum_{i=1}^m a_i\, c^{b, \theta}_i$ is defined by
$$
\partial \Big(\displaystyle\sum_{i=1}^m a_i\, c^{b, \theta}_i\Big)
:=\displaystyle\sum_{i=1}^m a_i \,\partial (c^{b, \theta}_i).
$$

Let $w$ be a dual real differentiable $k$-form on
$[0,\, b]_{\theta}^k\subseteq \mathcal{R}^{(2)k}$. Then there exists a unique dual real differentiable function $f(x^1, \dots , x^k)$ such that
$$w=f\, dx^1\wedge\cdots \wedge dx^k.$$
Using the uniqueness of the function $f(x^1, \dots , x^k)$ and the type $\theta$ integral introduced in section 3, we define
$$
\displaystyle\int _{[0,\, b]_{\theta}^k}\, w:=\displaystyle\int _{[0,\, b]_{\theta}^k}\,
f(x^1, \dots , x^k)\,d_{\theta}x^1\cdots d_{\theta}x^k,
$$
which is also written as
$$
\displaystyle\int _{[0,\, b]_{\theta}^k}\,f\, dx^1\wedge\cdots \wedge dx^k
=\displaystyle\int _{[0,\, b]_{\theta}^k}\,
f(x^1, \dots , x^k)\,d_{\theta}x^1\cdots d_{\theta}x^k .
$$

If $w$ is a dual real differentiable $k$-form on $A\subseteq \mathcal{R}^{(2)n}$ and
$c^{b, \theta}: [0,\, b]_{\theta}^k\to A$ is a type $\theta$ singular $k$-cube in $A$ with $k\ge 1$, we define
\begin{equation}\label{eq4.a}
\displaystyle\int _{c^{b, \theta}}\,w:=\displaystyle\int _{[0,\, b]_{\theta}^k}\,
(c^{b, \theta})^*\,w ,
\end{equation}
where $(c^{b, \theta})^*\,w $ is the dual real differentiable $k$-form on
$[0,\, b]_{\theta}^k$ defined by (\ref{eq4.70}).

\medskip
If $w$ is a  dual real differentiable $0$-form and $c^{b, \theta}: [0,\, b]_{\theta}^0=\{0\}\to A$ is (type $\theta$) singular $0$-cube in $A$, we define
$$
\displaystyle\int _{c^{b, \theta}} w:=w\Big(c^{b, \theta}(0)\Big) .
$$

The type $\theta$ integral of a dual real differentiable $k$-form over a type $\theta$ $k$-chain
$c^{b, \theta}=\displaystyle\sum_{i=1}^m a_i\, c^{b, \theta}_i$ is defined by
$$
\displaystyle\int _{c^{b, \theta}}w:=\displaystyle\sum_{i=1}^m a_i
\displaystyle\int _{c^{b, \theta}_i}w .
$$

\medskip
\begin{proposition}\label{pr4.13} ({\bf Type $\theta$ Stokes' Theorem for Type $\theta$ $k$-Chain}) If $w$ is a dual real differentiable $(k-1)$-form on an open set $A\subseteq \mathcal{R}^{(2)n}$ and $c^{b, \theta}$ is a type $\theta$ $k$-chain in $A$ with $k\ge 1$, then
$$
\displaystyle\int _{c^{b, \theta}}\, dw=\displaystyle\int _{\partial c^{b, \theta}}\, w .
$$
\end{proposition}

\end{document}